\documentclass{article}

\usepackage{PRIMEarxiv}

\usepackage[utf8]{inputenc} 
\usepackage[T1]{fontenc}    

\usepackage{amsmath}
\usepackage{amssymb}
\usepackage{bm}
\usepackage[numbers]{natbib}

\usepackage{hyperref}       
\usepackage{url}            
\usepackage{booktabs}       
\usepackage{amsfonts}       
\usepackage{nicefrac}       
\usepackage{microtype}      
\usepackage{graphicx}       
\usepackage{subfigure}      
\usepackage{multirow}
\usepackage{algorithm}
\usepackage{algpseudocode}
\graphicspath{{Figures/}{media/}} 

\title{Tiny Recursive Models for Solving the $J_2$-Perturbed Lambert Problem
\thanks{\textit{\underline{Citation}}:
\textbf{Wijayatunga, M, Armelin A. Tiny Recursive Models for Solving the $J_2$-Perturbed Lambert Problem}\\
Code: https://github.com/MinduliW/TRCLambert}
}

\author{
  Minduli Wijayatunga \\
  Department of Aerospace Engineering \\
  University of Illinois Urbana-Champaign \\
  Urbana-Champaign, IL, USA \\
  \texttt{email@email} \\
  \And
  Roberto Armellin \\
  Te P\=unaha \=Atea -- Space Institute \\
  University of Auckland \\
  Auckland, New Zealand \\
  \texttt{email@email} \\
}

\begin{document}
\maketitle

\begin{abstract}
This paper presents a fast, recursive neural solver for the $J_2$-perturbed Lambert problem based on Tiny Recursive Models (TRM), named TRM-Perturbed Lambert (TRM-PL) Model. TRM is a weight-shared neural architecture whose effective capacity emerges from iteration depth rather than parameter count. In TRM, a compact reasoning module is repeatedly applied within a two-level latent hierarchy, refining a candidate departure velocity by simulating the $J_2$ trajectory and correcting it based on the resulting tracking error.  This design unifies initial-guess generation and iterative correction in a single end-to-end differentiable architecture. 
The recursive refinement loop can be viewed as a learned, reinforcement-learning-like alternative to the homotopy and continuation schemes used in classical perturbed-Lambert solvers, where rather than following a hand-designed continuation path from the Keplerian to the perturbed solution, the network learns its own sequence of corrections. We evaluate TRM on four test cases of increasing difficulty: single- and multi-revolution low-Earth-orbit (LEO) transfers, and single- and multi-revolution Jovian transfers. Three training paradigms are compared: jointly learning the Lambert solution and the $J_2$ correction; refining the Lambert initial velocity with explicit target-position and $J_2$-corrected velocity supervision; and refining the Lambert initial velocity with only target-position supervision.  Across all four test cases, the refinement-only approaches emerge as the most reliable solver for the $J_2$ perturbed Lambert's problem. Across the three test cases, the position-supervised refinement variant reduces the median terminal-position error from 21.7\,km to 0.027\,km on single-revolution LEO, from 340.9\, km to 0.31\,km on multi-revolution LEO, and from $\sim$150\,km to 3.9\,km on multi-revolution Jovian, all using the same 2.3 M-parameter architecture. With one Newton corrector iteration applied to the TRM-PL output, the Jovian median tightens further to 0.063\,km. The resulting models are compact and accurate enough for embedded deployment, and can be refined further with Newton iterations as needed. 
\end{abstract}


\section{Introduction}
 
Lambert's problem entails determining the conic transfer arc connecting two position vectors for a prescribed time of flight. It is a foundational tool in mission design, orbit determination, space surveillance, and conjunction analysis~\cite{battin1999,bate1971}. Fast boundary-value solvers are also central to autonomous mission-design pipelines for multi-target rendezvous and debris-removal missions, where many candidate transfers must be evaluated under perturbations and operational constraints \cite{wijayatunga2023multiadr}. Under purely Keplerian dynamics, modern algorithms~\cite{izzo2015,armellin2018} solve it robustly and in microseconds for both single- and multi-revolution cases.   The dynamics of any real central body, however, departs from the Keplerian model. Among the various non-Keplerian effects, the $J_2$ zonal harmonic of the gravitational potential dominates for most planetary bodies, and on a single transfer arc its secular contribution is the term that most visibly breaks the Keplerian Lambert prediction. Even for Earth, $J_2$ produces secular drift in the right ascension of the ascending node and the argument of perigee that, accumulated over a single transfer that is one orbital period or longer, displaces the propagated terminal position by tens to hundreds of kilometers relative to the Keplerian prediction. For oblate planets such as Jupiter, where $J_2$ is more than an order of magnitude stronger than Earth's, the same effect grows to thousands of kilometers on multi-revolution arcs~\cite{yang2022}. In all of these settings, the Keplerian Lambert departure velocity is no longer the velocity that actually reaches the targeted position when the spacecraft is propagated under the true dynamics, and a perturbed solver must be used instead.

The classical approach to the perturbed Lambert problem treats it as a two-point boundary-value problem and applies a Newton-type shooting iteration till the target position is reached \cite{kraige1982,arora2015}.  Each iteration requires at least one full numerical propagation under the perturbed dynamics, alongside three more for the finite-difference Jacobian. Single-revolution near-circular transfers converge in a few iterations from the Keplerian guess, but the basin of attraction shrinks rapidly with the number of revolutions and the strength of the perturbation: long transfers exhibit highly sensitive end conditions and the Keplerian initial guess often falls outside the convergence region~\cite{woollands2015}. Homotopy methods~\cite{yang2015homotopy,armellin2018} mitigate this by gradually introducing the perturbation, but the number of homotopy steps grows with the number of revolutions. Higher-order expansion methods~\cite{alhulayil2018} accelerate the inner iterations but still require a sufficiently good initial guess. In fact, the quality of the initial velocity provided to the corrector remains the bottleneck in most existing solvers for the perturbed Lambert problem.

A natural way to provide better initial guesses for the shooting method is to learn them from data. Yang et al.~\cite{yang2022} trained a feedforward deep neural network (DNN) to predict the velocity correction $\Delta \mathbf{v}$ between the Keplerian and $J_2$-consistent solutions, using the DNN output to warm-start a finite-difference Newton corrector. Their work provides a careful study of input-output sample design, reporting that a sample form built around spherical-coordinate inputs and the residual $\Delta \mathbf{v}$ as output is markedly easier to learn than alternative parameterizations. With this sample form they achieved 100\% convergence on Jovian transfers from 0 to 10 revolutions and total computation times below 0.5\,s, significantly outperforming standard Newton shooting and homotopy. 
The architecture, however, has two structural limitations. First, the network produces a single feed-forward initial guess: it cannot improve its own output if more computation is available, and accuracy plateaus at the level set by the network's training error (in their results, an initial-guess velocity error of $\sim$5\,m/s and an initial-guess position error of $\sim$100\,km, before the shooting corrector closes the gap). Second, the network and the corrector remain separate components: the DNN is trained to reproduce $\Delta \mathbf{v}$ from offline shooting data and has no end-to-end visibility into the corrector that consumes its output, which limits how much of the convergence load can be moved onto the network.
 
Recent work in machine reasoning suggests a different approach. Tiny Recursive Models (TRM)~\cite{jolicoeurmartineau2025morerecursivereasoningtiny} demonstrate that, for problems whose natural solver is iterative, a compact network applied repeatedly under shared weights can outperform much larger feed-forward networks. TRMs have been used for solving nonlinear optimal control problems by Jain and Linares~\cite{jain2026trc}, where they showed that approximately $1.5$\,M parameters suffice to match near-optimal control costs on Van der Pol stabilization and powered descent. What makes TRM efficient is that a small reasoning module is repeatedly applied to refine an estimate, with each pass conditioned on the tracking error from the previous one. Capacity comes from iteration depth, not from parameter count, and the procedure plays the role of a Newton corrector with a learned, state-dependent inverse Jacobian.

The recursive refinement loop is, in effect, a reinforcement-learned continuation scheme. Classical perturbed-Lambert solvers reach the solution by following a continuation path \cite{yang2015homotopy,armellin2018} in which a sequence of intermediate problems bridges the gap between the Keplerian and the perturbed solution. Similar continuation and initialization ideas have also been used in indirect trajectory optimization, where  problem transformations and scaling constants improve convergence to fuel- or energy-optimal solutions~\cite{wijayatunga2023scaling}. The recursive solver replaces 
this hand-designed path with a learned one where each refinement iteration is one step  of a continuation sequence. The network learns that sequence by mimicking reinforcement-learning (RL): it produces a velocity correction, observes the terminal-position error returned by the differentiable propagator, and learns to reduce it, rather than being supervised against a fixed velocity label. 
When the network is additionally tasked with reproducing the Keplerian Lambert solution itself, as in Variant~A below, there is an added component of behavioral cloning before the refinement, where the network is trained to imitate the classical Lambert solver and used to warm-start the RL-like refinement loop.

This work adopts the TRM architecture to solve the $J_2$-perturbed Lambert problem. The Keplerian Lambert solver provides a high-quality but biased initial guess, and the gap to the true solution is exactly the residual that a refinement operator should learn to drive to zero. The terminal position error after a $J_2$ propagation gives a natural feedback signal, and the same correction rule applies at every iteration, just as a Newton update does. Cast in this form, a recursive neural solver replaces the DNN-plus-shooting pipeline of Yang et al.~\cite{yang2022} with a single end-to-end differentiable architecture in which the propagator, error feedback, and correction live in the same computation graph.
The contributions of this paper are threefold:

\begin{enumerate}
 \item \textbf{Unified recursive architecture for the perturbed Lambert problem.} The $J_2$-perturbed Lambert problem is formulated as iterative refinement around the Keplerian solution. The TRM-Perturbed Lambert Model (TRM-PL) instantiates the TRM two-head architecture for this problem.

  \item \textbf{Systematic comparison of three training paradigms.} Three variants of the same architecture are trained to learn how to best perform the $J_2$ refinement,
  jointly learning the Lambert solution and the $J_2$ correction; refining from the Keplerian Lambert solution with both position and velocity supervision; and refining from the Keplerian Lambert solution with position-only supervision.

  \item \textbf{Evaluation across three problem regimes.} The solver is tested on single and multi-revolution LEO transfers (10--90\,min, altitudes 300--2000\,km, up to 13 revolutions) and on multi-revolution Jovian transfers (TOF up to 10 orbital periods, perijove 5--30\,$R_J$, up to 11 revolutions, following \cite{yang2022}). 
\end{enumerate}
 
 
\section{Perturbed Lambert Problem Formulation}
\label{sec:problem}

\begin{figure}[hbt!]
    \centering
    \includegraphics[width=\linewidth]{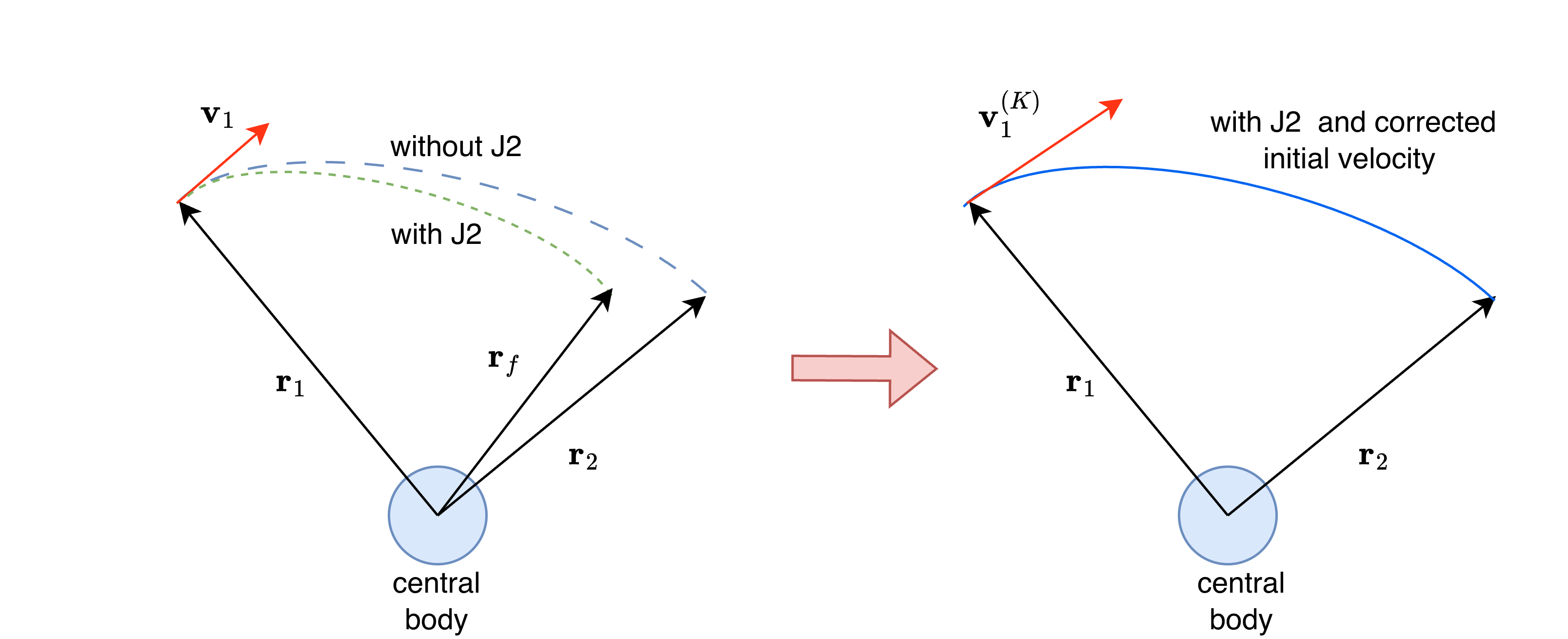}
    \caption{Illustration of the $J_2$-perturbed Lambert problem. The
    Keplerian Lambert arc (left, dashed green) connects $\mathbf{r}_1$
    to the target position $\mathbf{r}_2$ exactly under two-body dynamics;
    propagation under $J_2$ from the same initial condition
    $(\mathbf{r}_1,\mathbf{v}_{1,\mathrm{Lambert}})$ instead reaches the
    propagated terminal position $\mathbf{r}_f$, producing the terminal
    miss vector $\mathbf{r}_2-\mathbf{r}_f$ (left, dashed blue). The
    corrected initial velocity $\mathbf{v}_1^{(K)}$ produced by the
    recursive solver closes this gap.}
    \label{fig:TRMps}
\end{figure}

The perturbed Lambert problem considered in this work is illustrated in
Fig.~\ref{fig:TRMps}. Let $\mathbf{r}_1,\mathbf{r}_2 \in \mathbb{R}^3$
be the inertial-frame departure and desired arrival position vectors, and
let $\Delta t$ be the prescribed time of flight. A classical Lambert solver
under two-body dynamics returns a departure velocity
$\mathbf{v}_{1,\mathrm{Lambert}}$ such that the unperturbed trajectory
from $(\mathbf{r}_1,\mathbf{v}_{1,\mathrm{Lambert}})$ reaches
$\mathbf{r}_2$ exactly at time $\Delta t$.

Under $J_2$-perturbed dynamics, the equations of motion in an inertial
Earth-centered, or planet-centered, frame are
\begin{equation}
  \ddot{\mathbf{r}} = -\frac{\mu}{r^3}\mathbf{r} + \frac{3}{2}\frac{J_2 \mu R^2}{r^5}
  \begin{bmatrix}
    x\!\left(\dfrac{5z^2}{r^2} - 1\right) \\[4pt]
    y\!\left(\dfrac{5z^2}{r^2} - 1\right) \\[4pt]
    z\!\left(\dfrac{5z^2}{r^2} - 3\right)
  \end{bmatrix},
  \label{eq:j2eom}
\end{equation}
where $r = \|\mathbf{r}\|$, $\mu$ is the gravitational parameter, $R$ is
the equatorial radius of the central body, and $J_2$ is the leading zonal
harmonic coefficient.

Let $\mathbf{r}_f(\mathbf{r}_1,\mathbf{v}_1)$
denote the terminal position obtained by integrating Eq.~\eqref{eq:j2eom}
from the initial condition $(\mathbf{r}_1,\mathbf{v}_1)$ over the time of
flight $\Delta t$. In general, $ \mathbf{r}_f(\mathbf{r}_1,\mathbf{v}_{1,\mathrm{Lambert}})
    \neq
    \mathbf{r}_2,$
because $\mathbf{v}_{1,\mathrm{Lambert}}$ is computed under two-body
dynamics, whereas the propagated trajectory includes the $J_2$ perturbation.
The corresponding terminal-position error is defined as
\begin{equation}
    \mathbf{e}(\mathbf{v}_1)
    =
    \mathbf{r}_2
    -
    \mathbf{r}_f(\mathbf{r}_1,\mathbf{v}_1).
    \label{eq:terminal_error}
\end{equation}
Thus, the initial Lambert miss under $J_2$ propagation is
\begin{equation}
    \mathbf{e}_{\mathrm{Lambert}}
    =
    \mathbf{r}_2
    -
    \mathbf{r}_f(\mathbf{r}_1,\mathbf{v}_{1,\mathrm{Lambert}}).
\end{equation}

The $J_2$-perturbed Lambert problem  seeks a velocity correction
$\Delta\mathbf{v}$ such that
\begin{equation}
  \big\|
  \mathbf{r}_2
  -
  \mathbf{r}_f(\mathbf{r}_1,
  \mathbf{v}_{1,\mathrm{Lambert}} + \Delta\mathbf{v})
  \big\|
  \leq \varepsilon_r,
  \label{eq:j2lambert}
\end{equation}
for a prescribed terminal-position tolerance $\varepsilon_r$. The corrected
solution is denoted
\begin{equation}
    \mathbf{v}_{1,\mathrm{true}}
    =
    \mathbf{v}_{1,\mathrm{Lambert}} + \Delta\mathbf{v}.
\end{equation}

 \section{Newton Shooting Method Baseline}\label{ns}

The classical approach to the perturbed Lambert problem treats it as a
two-point boundary-value problem and applies a Newton-type shooting
iteration until the propagated terminal position reaches the target
position. This method is given in Algorithm~\ref{alg:newton_j2_shooting}.
It computes the $J_2$-corrected departure velocity by linearizing the
terminal propagation map around the current velocity iterate. In this work, the performance of TRM-PL is compared against this baseline.


\begin{algorithm}[hbt!]
\caption{Newton shooting for the $J_2$-corrected departure velocity}
\label{alg:newton_j2_shooting}
\begin{algorithmic}[1]
\Require Initial position $\mathbf{r}_1$, target position $\mathbf{r}_2$, time of flight $\Delta t$, initial guess $\mathbf{v}_{1,\mathrm{Lambert}}$, tolerance $\varepsilon_r$, finite-difference step $\epsilon_v$,  maximum iterations $j_{\max}$
\Ensure Corrected departure velocity $\mathbf{v}_{1,\mathrm{true}}$
\State Initialise $\mathbf{v}_1^{(0)} \gets \mathbf{v}_{1,\mathrm{Lambert}}$
\For{$j = 0,\dots,j_{\max}-1$}
    \State Propagate Eq.~\eqref{eq:j2eom} from $(\mathbf{r}_1,\mathbf{v}_1^{(j)})$ over $\Delta t$
    \State Compute the propagated terminal position $\mathbf{r}_f^{(j)} = \mathbf{r}_f(\mathbf{r}_1,\mathbf{v}_1^{(j)})$
    \State Compute the terminal-position error as  $\mathbf{e}^{(j)}
        = \mathbf{r}_2  - \mathbf{r}_f^{(j)}$. 
        
    \If{$\|\mathbf{e}^{(j)}\| \leq \varepsilon_r$}
        \State \Return $\mathbf{v}_1^{(j)}$
    \EndIf
    \State Approximate the error Jacobian by forward differences:
    \[
        \mathbf{S}^{(j)}_{:,k}
        \approx
        \frac{
        \mathbf{e}(\mathbf{v}_1^{(j)}+\epsilon_v\mathbf{q}_k)
        -
        \mathbf{e}(\mathbf{v}_1^{(j)})
        }{\epsilon_v},
        \qquad k=1,2,3
    \]
    where $\mathbf{q}_k$ is the $k$th Cartesian unit vector in velocity space.
    \State Solve the Newton system $\Delta\mathbf{v}^{(j)}
        =
        -\left(\mathbf{S}^{(j)}\right)^{-1}\mathbf{e}^{(j)}$.
    \State Update the departure velocity: $\mathbf{v}_1^{(j+1)}
        =
        \mathbf{v}_1^{(j)}
        +
        \Delta\mathbf{v}^{(j)}$
\EndFor
\State \Return failure
\end{algorithmic}
\end{algorithm}

 \section{Tiny Recursive Models for the J2-Perturbed Lambert Problem}

The structure of the terminal-position error, $\mathbf{e}$, motivates the use of a TRM for the perturbed Lambert problem. This error is a nonlinear function of the boundary data $(\mathbf{r}_1,\mathbf{r}_2,\Delta t)$ and the initial Lambert velocity $\mathbf{v}_{1,\mathrm{Lambert}}$. Classical Newton shooting reduces this error by repeatedly approximating the local map from changes in initial
velocity to changes in terminal position; however, it requires at least four numerical propagations under perturbed dynamics per iteration. The proposed TRM instead learns the inverse refinement map from  $\mathbf{e}$ to $\Delta\mathbf{v}$ directly by predicting the velocity correction that moves $\mathbf{e}$ to zero given the current error and problem context. The learned
recursive update plays the role of a data-driven Newton correction without
explicitly computing or inverting the sensitivity matrix.

\subsection{TRM Solver Architecture}

\begin{figure}[hbt!]
    \centering
    \subfigure[Variant A: learned Lambert.]{%
        \includegraphics[width=0.45\linewidth]{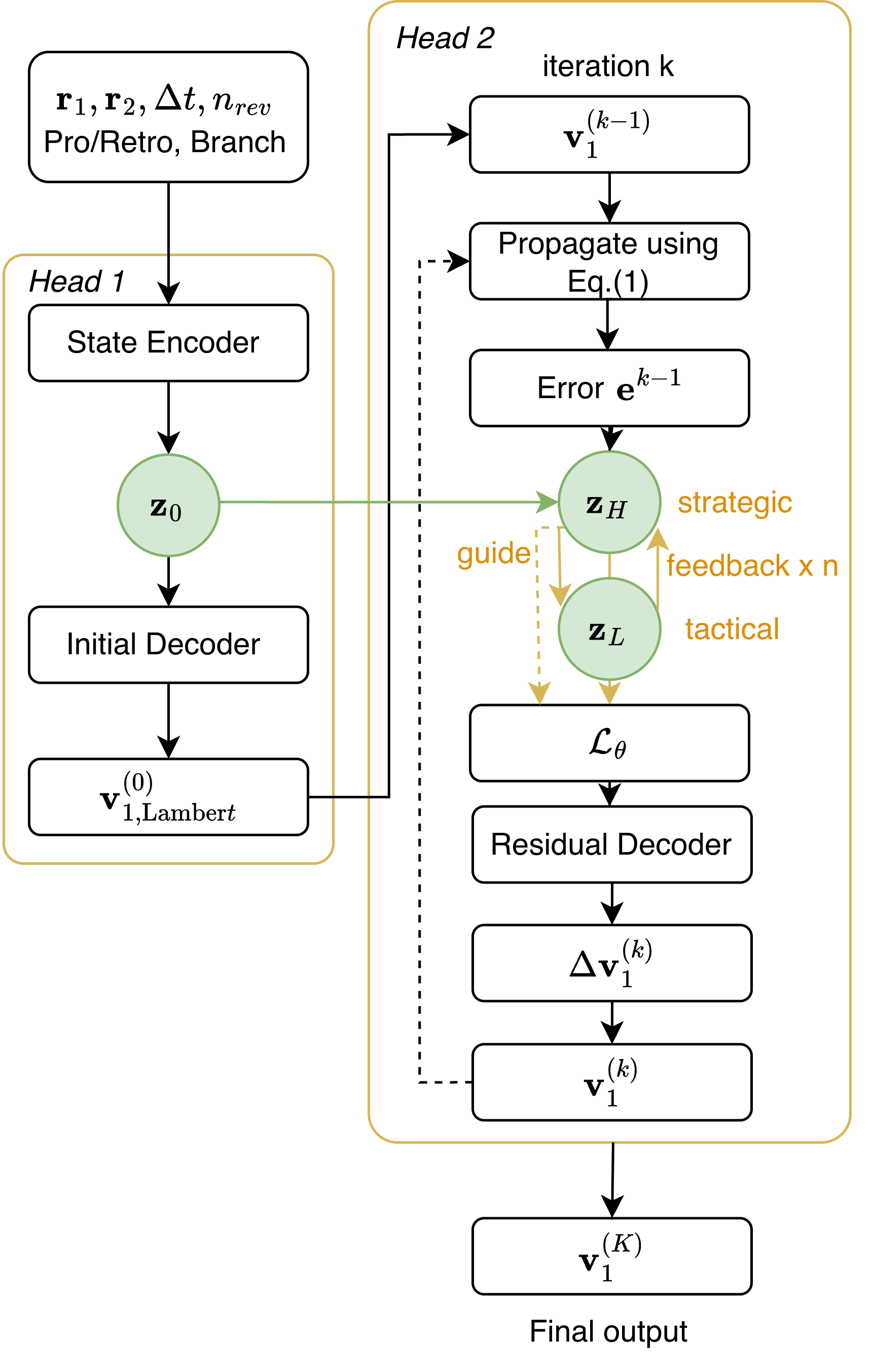}%
        \label{fig:TRMmethod_a}%
    }\hfill
    \subfigure[Variants B and C: refinement-only.]{%
        \includegraphics[width=0.45\linewidth]{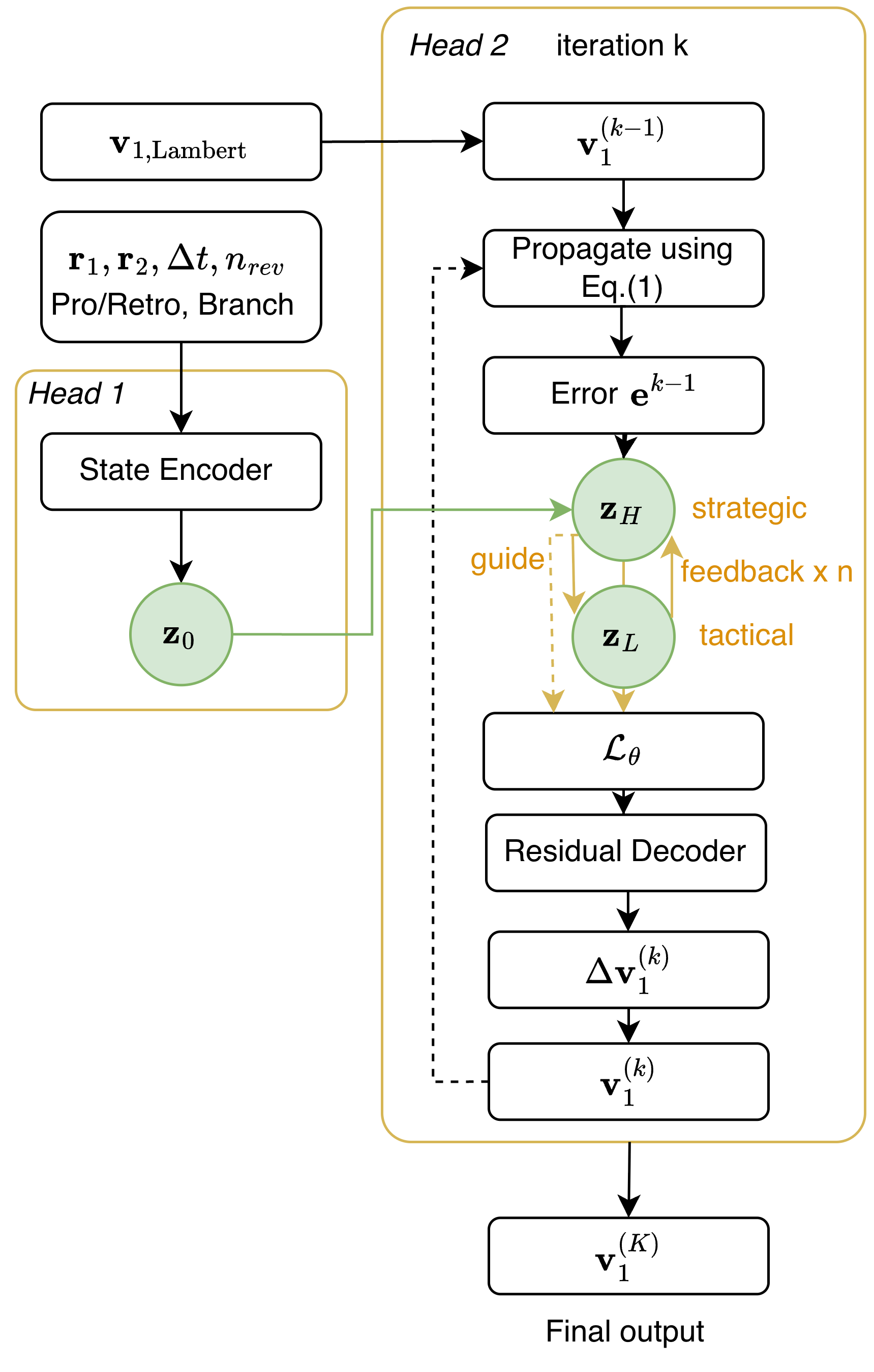}%
        \label{fig:TRMmethod_bc}%
    }
    \caption{Two-head TRM architecture for the $J_2$-perturbed Lambert problem. The two configurations differ only in the source of $\mathbf{v}_1^{(0)}$: a learned Initial Decoder in (a), the classical Lambert solver in (b).}
    \label{fig:TRMmethod}
\end{figure}

TRM-PL, illustrated in Figure~\ref{fig:TRMmethod}, is hereafter referred to as the TRM-Perturbed Lambert Model (TRM-PL). It instantiates the two-head TRM architecture of~\citep{jolicoeurmartineau2025morerecursivereasoningtiny} for the $J_2$-perturbed Lambert problem. Inputs are the boundary conditions $(\mathbf{r}_1, \mathbf{r}_2, \Delta t)$ together with the discrete Lambert case parameters $(n_\mathrm{rev}, \mathrm{prograde}, \mathrm{branch})$; the output is the $J_2$-corrected departure velocity $\mathbf{v}_1^{(K)}$. Three variants of TRM-PL are explored to determine the most efficient solution method. They are 
 \begin{itemize}
 \item \textbf{Variant A -- Learned Lambert.}
Head 1 is supervised against the Keplerian Lambert solution and is responsible for producing the initial guess $\mathbf{v}_1^{(0)}$ that Head 2 refines; Head 2 is supervised against the $J_2$-corrected solution. At inference time, the classical Lambert solver is not invoked at all, and a single unified network handles both initial-guess generation and iterative correction. In this configuration Head 1 acts as a behavioral-cloning warm start that initializes the reinforcement-learning-like refinement performed by Head 2. 


     \item \textbf{Variant B -- Velocity and Position-Supervised Refinement} 

     This variant uses the Lambert velocity as input to the refinement process, rather than using head 1 to guess the initial value of $\mathbf{v}_1^{(0)}$. The training is done by supervising both the terminal position error and the corrected velocity error. 

      \item \textbf{Variant C -- Position-Supervised Refinement}

      This Variant is similar to Variant B, however the training is done by supervising only using the terminal position error. No $J_2$ corrected velocities are required in the training data.

 \end{itemize}

 In all variants, then Head~2 applies $K$ recursive correction steps to produce $\mathbf{v}_1^{(K)}$, using a shared reasoning module that operates on a high/low-level latent hierarchy. The remainder of this subsection describes each component in turn: the state encoder, Head~1, the differentiable $J_2$ propagator that supplies the tracking error to Head~2, and the recursive refinement loop of Head~2 itself.
 
\subsubsection{State Encoder}

The state encoder $\mathrm{MLP}_{\mathrm{state}}$ embeds the boundary-value problem into a $d_z$-dimensional latent vector
\begin{equation}
    \mathbf{z}_0 = \mathrm{MLP}_{\mathrm{state}}\bigl([\mathbf{r}_1; \mathbf{r}_2; \Delta t; \mathrm{prograde}; n_\mathrm{rev}; \mathrm{branch}]\bigr),
    \label{eq:z0}
\end{equation}
where $\mathrm{MLP}_{\mathrm{state}}$ is a two-layer perceptron with the structure $\mathrm{Linear}\to\mathrm{LayerNorm}\to\mathrm{GELU}\to\mathrm{Linear}$. Position vectors are normalized by a per-dataset characteristic length and time of flight by a characteristic time before being passed to the encoder. The same encoding $\mathbf{z}_0$ is reused as a persistent problem-context token at every refinement iteration of Head~2, and feeds the Initial Decoder of Head~1 in Variant~A.

\subsubsection{Initial-Guess Generation (Variant A only)}

Head~1 produces a first estimate of the departure velocity directly from $\mathbf{z}_0$:
\begin{equation}
    \mathbf{v}_1^{(0)} = \mathrm{MLP}_{\mathrm{initial}}(\mathbf{z}_0),
    \label{eq:v0}
\end{equation}
where $\mathrm{MLP}_{\mathrm{initial}}$ has the same two-layer structure as the state encoder, with output dimension three. Predictions are made in normalized space and de-normalized by a velocity scale before being used. Head~1 plays the role of an explicit learned Lambert solver, a feedforward replacement for the classical algorithm of Izzo~\cite{izzo2015}. 


\subsubsection{J2 Propagation and Error Computation}
 
A differentiable batch propagator is used to implement Eq.~\eqref{eq:j2eom} via classical fourth-order Runge--Kutta integration. The number of integration steps is set per sample as $n_\mathrm{step} = \max(50,\lceil \Delta t / h_\mathrm{max}\rceil)$, where $h_\mathrm{max}$ is a problem-dependent maximum step size.  At iteration $k$, the propagator returns $\mathbf{r}_f^{(k-1)},\, \mathbf{v}_f^{(k-1)}$,
from which the terminal-position residual is computed as
\begin{equation}
     \mathbf{e}^{(k-1)} = \mathbf{r}_f^{(k-1)} - \mathbf{r}_2.
     \label{eq:err}
\end{equation}

\subsubsection{Head 2: Iterative Refinement}
 
The refinement head consumes either the estimated Lambert velocity $\mathbf{v}_1^{(0)}$ for Variant A or the true Lambert velocity $\mathbf{v}_{1,\mathrm{Lambert}}$ for Variants B and C, 
and applies $K$ correction iterations to produce successively more accurate estimates $\mathbf{v}_1^{(1)}, \ldots, \mathbf{v}_1^{(K)}$. Two latent states are carried across iterations: a high-level (\emph{strategic}) latent $\mathbf{z}_H \in \mathbb{R}^{d_z}$ and a low-level (\emph{tactical}) latent $\mathbf{z}_L \in \mathbb{R}^{d_z}$. Following~\cite{jolicoeurmartineau2025morerecursivereasoningtiny}, both are initialized as a sum of a globally learnable parameter and a sample-specific projection of the problem encoding:
\begin{equation}
    \mathbf{z}_H^{(0)} = \mathbf{H}_{\mathrm{init}} + \mathbf{W}_H \mathbf{z}_0,
    \qquad
    \mathbf{z}_L^{(0)} = \mathbf{L}_{\mathrm{init}} + \mathbf{W}_L \mathbf{z}_0,
\end{equation}
where $\mathbf{H}_{\mathrm{init}},\mathbf{L}_{\mathrm{init}} \in \mathbb{R}^{d_z}$ are global learnable vectors and $\mathbf{W}_H,\mathbf{W}_L$ are learned linear projections.  The same reasoning weights are reused at every iteration,leading to the efficiency of TRM.
At each outer iteration $k = 1, \ldots, K$, the position residual $\mathbf{e}^{(k-1)}$ from Eq.~\eqref{eq:err} and the current velocity iterate $\mathbf{v}_1^{(k-1)}$ are encoded into latent tokens
\begin{equation}
\mathbf{z}_{\mathrm{err}}^{(k)} = \mathrm{MLP}_{\mathrm{err}}\bigl(\mathbf{e}^{(k-1)} \bigr),
\qquad
\mathbf{z}_{\mathrm{ctrl}}^{(k)} = \mathrm{Linear}\bigl(\mathbf{v}_1^{(k-1)}\bigr),
\label{eq:errctrl}
\end{equation}
 The low-level latent then performs $n$ tactical cycles, each updating $\mathbf{z}_L$ from the full token set
\begin{equation}
  \mathbf{z}_L^{(i+1)} = \mathcal{L}_\theta\!\bigl(
    \mathbf{z}_L^{(i)},\; \mathbf{z}_H^{(k-1)},\; \mathbf{z}_0,\;
    \mathbf{z}_{\mathrm{err}}^{(k)},\; \mathbf{z}_{\mathrm{ctrl}}^{(k)}\bigr),
    \qquad i = 0, \ldots, n-1,
    \label{eq:zL}
\end{equation}
where the first argument is the updated output token. This formulation differs from the original Tiny Recurrent Control (TRC) paper~\cite{jain2026trc}, which sums the auxiliary context into a single vector before passing it to the reasoning module. Preserving the five tokens individually allows the multi-head self-attention layers to weight problem context, strategic memory, error feedback, and current control separately, which was found to consistently lower the validation error in the Lambert setting.

After the $n$ tactical cycles complete, the high-level latent updates once, based on the new tactical state as 
\begin{equation}
    \mathbf{z}_H^{(k)} = \mathcal{L}_\theta\bigl( \mathbf{z}_H^{(k-1)}, \mathbf{z}_L^{(n)}\bigr).
    \label{eq:zH}
\end{equation}
 
The shared reasoning module $\mathcal{L}_\theta$ is a stack of $L = 3$ transformer-style blocks, each consisting of multi-head self-attention with $h = 8$ heads followed by a feedforward sublayer, with residual connections and pre-LayerNorm normalization. The same weights are used in Eqs.~\eqref{eq:zL} and~\eqref{eq:zH}: the network learns a single context-conditioned update rule rather than separate strategic and tactical operators, and the role each application plays is determined by which tokens are present at its input.
 
Finally, the residual decoder produces a velocity correction conditioned on the updated strategic latent and the current velocity as
\begin{equation}
    \Delta\mathbf{v}^{(k)} = \mathrm{MLP}_{\mathrm{residual}}\bigl([\mathbf{z}_H^{(k)};\, \mathbf{v}_1^{(k-1)}]\bigr),
    \label{eq:dv}
\end{equation}
and the velocity is updated by residual addition with clipping to enforce a velocity bound,
\begin{equation}
\mathbf{v}_1^{(k)} = \mathrm{clip}\bigl(\mathbf{v}_1^{(k-1)} + \Delta\mathbf{v}^{(k)},\, -v_\mathrm{max},\, v_\mathrm{max}\bigr),
\label{eq:vupdate}
\end{equation}
 After the loop, a final $J_2$ propagation is performed with $\mathbf{v}_1^{(K)}$ to compute the terminal position used in the loss; this aligns the position term and velocity term of the loss to the same iterate.
 
The complete forward pass is summarized in Algorithm~\ref{alg:trm_lambert}.
 
\begin{algorithm}[t]
\caption{TRM-PL recursive solver (Head 1 + Head 2)}
\label{alg:trm_lambert}
\begin{algorithmic}[1]
\Require Boundary conditions $(\mathbf{r}_1,\mathbf{r}_2,\Delta t)$ and Lambert case $(n_\mathrm{rev},\mathrm{prograde},\mathrm{branch})$; trained TRM-PL weights $\theta$; iteration count $K$; optional oracle velocity $\mathbf{v}_1^{\mathrm{oracle}}$
\Ensure Refined departure velocity $\mathbf{v}_1^{(K)}$
\State $\mathbf{z}_0 \gets \mathrm{MLP}_{\mathrm{state}}([\mathbf{r}_1;\mathbf{r}_2;\Delta t;\mathrm{prograde};n_\mathrm{rev};\mathrm{branch}])$ \Comment{problem encoding, Eq.~\eqref{eq:z0}}
\If{Variant~A}
    \State $\mathbf{v}_1^{(0)} \gets \mathrm{MLP}_{\mathrm{initial}}(\mathbf{z}_0)$ \Comment{Initial Decoder, Eq.~\eqref{eq:v0}}
\Else \Comment{Variants~B and~C}
    \State $\mathbf{v}_1^{(0)} \gets \mathbf{v}_{1,\mathrm{Lambert}}$ \Comment{Classical Lambert solver}
\EndIf
\State $\mathbf{z}_H \gets \mathbf{H}_\mathrm{init} + \mathbf{W}_H \mathbf{z}_0$;\quad $\mathbf{z}_L \gets \mathbf{L}_\mathrm{init} + \mathbf{W}_L \mathbf{z}_0$
\For{$k = 1, \ldots, K$}
    \State $(\mathbf{r}_f^{(k-1)},\mathbf{v}_f^{(k-1)}) \gets \Phi_{\Delta t}^{J_2}(\mathbf{r}_1,\mathbf{v}_1^{(k-1)})$ \Comment{differentiable RK4}
    \State $\mathbf{e}^{(k-1)} \gets \mathbf{r}_f^{(k-1)} - \mathbf{r}_2$
    \State $\mathbf{z}_\mathrm{err} \gets \mathrm{MLP}_\mathrm{err}(\mathbf{e}^{(k-1)}/\sigma_r)$;\quad $\mathbf{z}_\mathrm{ctrl} \gets \mathrm{Linear}(\mathbf{v}_1^{(k-1)})$
    \For{$i = 1, \ldots, n$} \Comment{tactical cycles}
        \State $\mathbf{z}_L \gets \mathcal{L}_\theta(\mathbf{z}_L, \mathbf{z}_H, \mathbf{z}_0, \mathbf{z}_\mathrm{err}, \mathbf{z}_\mathrm{ctrl})$
    \EndFor
    \State $\mathbf{z}_H \gets \mathcal{L}_\theta(\mathbf{z}_H, \mathbf{z}_L)$ \Comment{strategic update}
    \State $\Delta\mathbf{v}^{(k)} \gets \mathrm{MLP}_\mathrm{residual}([\mathbf{z}_H; \mathbf{v}_1^{(k-1)}])$
    \State $\mathbf{v}_1^{(k)} \gets \mathrm{clip}(\mathbf{v}_1^{(k-1)} + \Delta\mathbf{v}^{(k)}, -v_\mathrm{max}, v_\mathrm{max})$
\EndFor
\State \Return $\mathbf{v}_1^{(K)}$
\end{algorithmic}
\end{algorithm}
 
\subsection{Training Methodology}

\subsubsection{Loss Functions}
\label{sec:variants}
The three variants share the same architecture and the same propagator. They differ only in (i) whether the Initial Decoder is trained or replaced by the classical Lambert solver, and (ii) whether the corrected departure velocity is supervised.

For Variant A,  Head 1 is trained against the Keplerian Lambert solution with the loss
\begin{equation}
    \mathcal{L}_\mathrm{Head \,1}(A) = \big\| \mathbf{v}_{1,\mathrm{Lambert}} - \mathbf{v}_1^{(0)} \big\|^2,
    \label{eq:stage1loss}
\end{equation}

Then, for variants A and B,  Head~2 loss is defined as follows:
\begin{equation}
\label{s2rew}
    \mathcal{L}_\mathrm{Head\,2}(A/B) \;=\;
\lambda_v\, \big\|\mathbf{v}_1^{(K)} - \mathbf{v}_{1,\mathrm{true}}\big\|^2
\;+\; \lambda_\mathrm{pos}\, \big\|\mathbf{e}^{(K)}\big\|^2
\;-\; \lambda_\mathrm{ps}\, \frac{1}{K-1}\sum_{k=1}^{K-1}\bigl(\tilde J^{(k-1)} - \tilde J^{(k)}\bigr),
\end{equation}
penalizing the deviation from the Newton-shooting velocity label, the final terminal-position error, and any iteration that fails to reduce the normalized terminal cost, respectively. For variant C, the Head~2 loss is
\begin{equation}
\label{s2rew2}
    \mathcal{L}_\mathrm{Head\,2}(C) \;=\; \lambda_\mathrm{pos}\, \big\|\mathbf{e}^{(K)}\big\|^2
\;-\; \lambda_\mathrm{ps}\, \frac{1}{K-1}\sum_{k=1}^{K-1}\bigl(\tilde J^{(k-1)} - \tilde J^{(k)}\bigr).
\end{equation}

 Note that 
\begin{equation}
J^{(k)} \;=\; \tfrac{1}{2}\,\bigl\|\mathbf{e}^{(k)}\bigr\|^2,
\qquad
\tilde J^{(k)} \;=\; \frac{J^{(k)}}{J^{(0)}}.
\label{eq:Jk}
\end{equation}
 The negative sign on the third term in Eq.turns it into a reward for monotonic decrease of $\tilde J$ across iterations. All three variants use $\lambda_\mathrm{pos} = 1.0$ and $\lambda_\mathrm{ps} = 0.1$. Variants~A and~B use $\lambda_v = 1.0$.

 Variant~C uses $\lambda_v = 0.0$, removing the velocity-label term entirely.

Training uses the AdamW optimizer with cosine learning-rate annealing, gradient clipping at norm 1.0, and batch size 512. 
Variants B and C train Head~2 alone for 300 epochs from a fresh initialization at learning rate $10^{-4}$; Variant A trains Head 1 to convergence first (500 epochs at learning rate $10^{-3}$ annealed to $10^{-6}$) and then trains Head 2 for an additional 300 epochs at $10^{-5}$. {Total wall-clock time is approximately 6--8 hours per run on a single NVIDIA RTX 3080 GPU for the LEO datasets and 12--16 hours for the Jovian datasets, dominated by RK4 propagation through the long Jovian time horizons.}

\subsubsection{Training Data Generation}
\label{sec:datagen}
All training and validation datasets are generated using the forward-sampling procedure of Yang et al.~\cite{yang2022}. Orbital elements $(a, e, i, \Omega, \omega, M_0)$ are drawn uniformly from case-specific ranges given in Table~\ref{samplingrange} and converted to the Cartesian initial state $(\mathbf{r}_1, \mathbf{v}_{1,\mathrm{true}})$. The time of flight $\Delta t$ is drawn uniformly from a case-specific window. The state is propagated under $J_2$ dynamics with a fixed-step RK4 integrator to obtain the terminal state $(\mathbf{r}_2, \mathbf{v}_{2,\mathrm{true}})$. The classical Keplerian Lambert problem is then solved on $(\mathbf{r}_1, \mathbf{r}_2, \Delta t)$ to obtain $\mathbf{v}_{1,\mathrm{Lambert}}$, with the appropriate Lambert branch identified directly from the propagated trajectory: $n_\mathrm{rev}$ is read off the trajectory, the prograde/retrograde flag is determined by the sign of the $z$-component of the orbital angular momentum $\mathbf{r}_1 \times \mathbf{v}_{1,\mathrm{true}}$, and the short- or long-way branch is selected as the Lambert solution that minimizes $\|\mathbf{v}_{1,\mathrm{Lambert}} - \mathbf{v}_{1,\mathrm{true}}\|_2$. Samples are rejected if either trajectory intersects the central body. Each retained sample stores the record $(\mathbf{r}_1, \mathbf{r}_2, \Delta t, n_\mathrm{rev}, \mathrm{prograde}, \mathrm{branch}, \mathbf{v}_{1,\mathrm{Lambert}}, \mathbf{v}_{1,\mathrm{true}})$.

\begin{table}[tbp]
\centering
\caption{Sampling ranges used to generate the three datasets, following the parameterization of Yang et al.~\cite{yang2022}. Semi-major axis and eccentricity are derived from the perijove/perigee and apojove/apogee radii.}
\label{samplingrange}
\small
\begin{tabular}{lccc}
\toprule
Variable & Single-rev LEO & Multi-rev LEO & Multi-rev Jovian \\
\midrule
Perigee/perijove radius $r_p$ & $[R_E + 300,\, R_E + 2000]$\,km & $[R_E + 300,\, R_E + 2000]$\,km & $[5,\,30]\,R_J$ \\
Apogee/apojove radius $r_a$   & $[r_p,\, R_E + 2000]$\,km        & $[r_p,\, R_E + 2000]$\,km        & $[r_p,\, 30]\,R_J$ \\
Inclination $i$ (rad)         & $[0,\,\pi]$                       & $[0,\,\pi]$                       & $[0,\,1]$ \\
RAAN $\Omega$ (rad)           & $[0,\,2\pi)$                      & $[0,\,2\pi)$                      & $[0,\,2\pi)$ \\
Argument of perigee $\omega$ (rad) & $[0,\,2\pi)$                 & $[0,\,2\pi)$                      & $[0,\,2\pi)$ \\
Mean anomaly $M_0$ (rad)      & $[0,\,2\pi)$                      & $[0,\,2\pi)$                      & $[0,\,2\pi)$ \\
Time of flight $\Delta t$ (T)    & (0,1)                  & (0,10)            & (0,10)  \\
\bottomrule
\end{tabular}
\end{table}

\section{Results}

This section presents the results obtained for the three test cases. The physical constants for LEO cases are  $\mu_E = 398{,}600.4418$\,km$^3$/s$^2$, $R_E = 6378.137$\,km, and $J_{2,E} = 1.08263 \times 10^{-3}$. For Jupiter, they are $\mu_J = 1.26686534\times 10^{8}$\,km$^3$/s$^2$, $R_J = 71{,}492$\,km, and $J_{2,J} = 1.4736\times 10^{-2}$. The differentiable RK4 integrator uses a maximum step size $h_\mathrm{max}$ chosen per dataset to balance accuracy against memory: $30$\,s for single-revolution LEO, $45$\,s for multi-revolution LEO, and $3600$\,s for Jovian transfers. The actual step count per sample is $n_\mathrm{step} = \max(50, \lceil \Delta t / h_\mathrm{max} \rceil)$.


All three test cases use the same network with latent dimension $d_z = 256$, hidden dimension $d_h = 512$, $L = 3$ reasoning blocks with $h = 8$ attention heads, and $n = 6$ inner cycles per outer iteration. The number of outer iterations $K$ is fixed per problem class at evaluation time, at $K = 3$ for single-revolution LEO and $K = 4$ for the multi-revolution LEO and multi-revolution Jovian cases. For all test cases, Variant~A has 2.6\,M trainable parameters; Variants~B and~C are smaller at 2.3\,M as the Initial Decoder is omitted. By comparison, the feedforward DNN of Yang et al.~\cite{yang2022} contains four hidden layers of 50 neurons each ($\sim$10\,K parameters) but produces only a single feedforward initial guess and relies on an external Newton corrector. Modern transformer-based control approaches~\cite{celestini2024,brohan2023rt2} use 50\,M to several billion parameters for problems of comparable structural complexity.

 \subsection{Training and Validation Data Distributions}
The mean, standard deviation, and magnitude difference coefficients of the time of flight, $\mathbf{v}_{1,\mathrm{Lambert}}, $ $\mathbf{v}_{1,\mathrm{true}}$ and $\mathbf{e}_{\mathrm{Lambert}}$ of the three resulting datasets are summarized in Appendix Table~\ref{tab:stats_all} and \ref{tab:stats_all_val}.

\subsection{Outcomes}
This section evaluates the three training variants on the held-out validation sets of all three test cases. Throughout, terminal-position error is reported as $\|\mathbf{r}_2 - \mathbf{r}_f^{(K)}\|$ (the distance between the propagated terminal position at the final iteration and the target) and corrected initial-velocity error as $\|\mathbf{v}_{1,\mathrm{true}} - \mathbf{v}_1^{(K)}\|$. Per-component statistics are also tabulated. All variants use the same architecture, the same propagator, and the same iteration count $K$ at evaluation time.

\subsubsection{Single Revolution LEO Results}

\begin{figure}[tbp]
    \centering
    \includegraphics[width=1\linewidth]{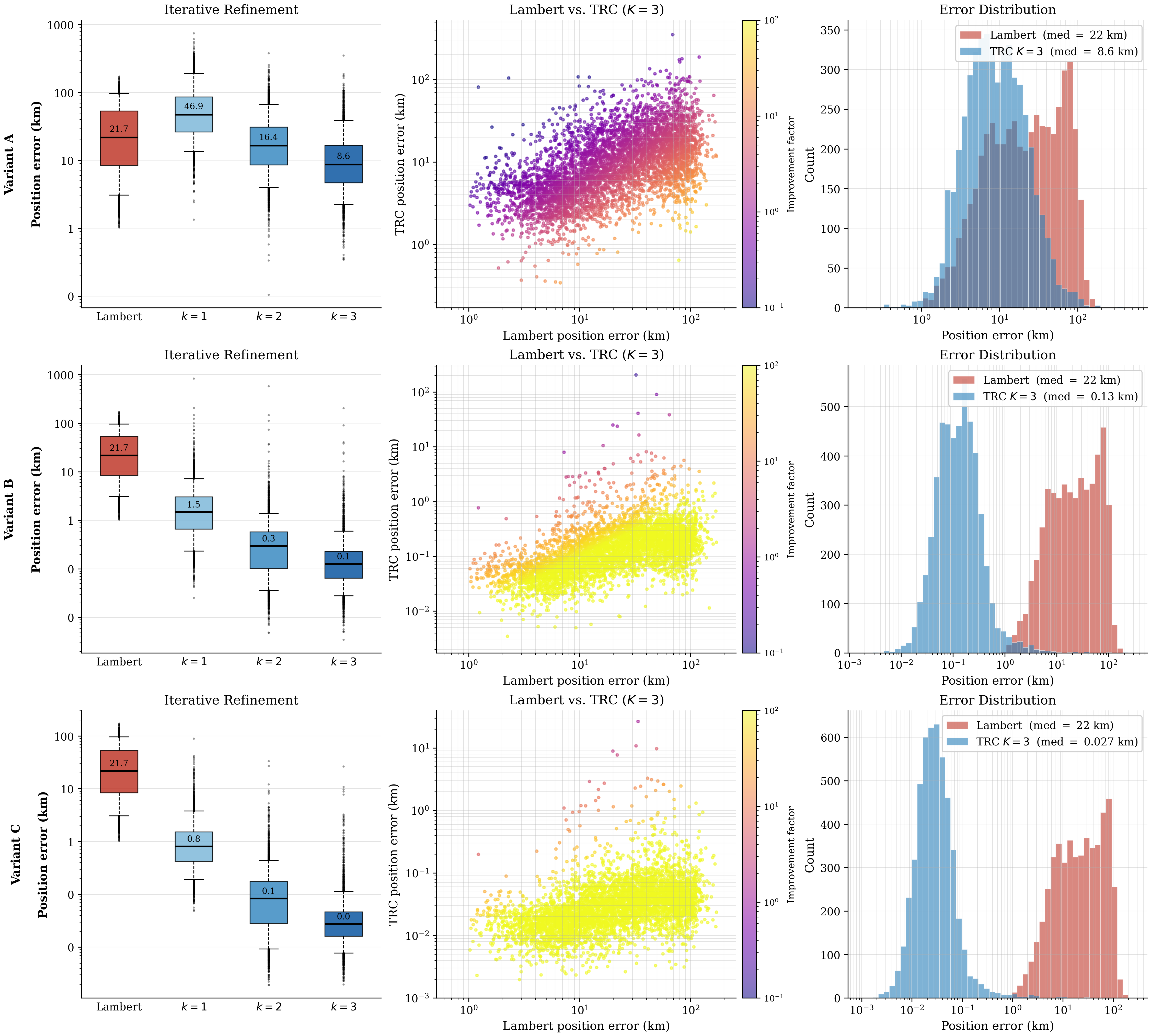}
        \caption{Performance of TRM on the single-revolution LEO dataset}
    \label{fig:s2_singlerev}
\end{figure}

The three variants are evaluated on the single-revolution LEO validation set of $N = 5$k $J_2$-perturbed Lambert transfers. Results are reported in Table~\ref{tab:err_single_rev} and Figure~\ref{fig:s2_singlerev}. All three variants reduce the terminal error relative to the uncorrected Keplerian Lambert baseline, but the refinement-based variants are more accurate by two orders of magnitude. Variant~C (position-supervised refinement) reaches a median terminal-position error of 0.027\,km and a mean of 0.065\,km, with the 99th percentile remaining below 0.52\,km. Variant~B (position- and velocity-supervised refinement) is comparable, with a median of 0.12\,km and a mean of 0.234\,km. The end-to-end learned-Lambert variant (Variant~A) is substantially less accurate, with a median of 16.7\,km and a mean of 31.8\,km, worse than Variant~B by roughly two orders of magnitude.

The same ordering holds in corrected-velocity error. Variant~C reaches a median velocity error of 0.020\,m/s, Variant~B 0.087\,m/s, and Variant~A 13.8\,m/s. The gap between Variants~B and~C is small in the central tendency but widens in the tail: Variant~C's 99th-percentile velocity error is 9.27\,m/s versus 25.9\,m/s for Variant~B, suggesting that the velocity-supervision term in Variant~B's loss adds variance without improving the upper tail of the distribution.

\begin{table}[tbp]
\centering
\small
\setlength{\tabcolsep}{5pt}
\caption{Final position and corrected initial velocity error distributions on the validation set of single-revolution LEO. $Q_2$ is the median; $P_{99}$ is the 99th percentile.}
\label{tab:err_single_rev}
\begin{tabular}{lrrrrr}
\toprule
Quantity & Mean & $Q_1$ & $Q_2$ & $Q_3$ & $P_{99}$ \\
\midrule
\multicolumn{6}{l}{\textbf{Variant A - Learned Lambert}}  \\
$\parallel \mathbf{r}_2 - \mathbf{r}_f \parallel$ (km) & 13.5 & 4.68 & 8.62 & 16.8 & 77.6 \\
$\parallel \mathbf{v}_{1,true} - \mathbf{v}_{1}^{(K)} \parallel$ (m/s) & 15.7 & 3.11 & 5.57 & 11.3 & 163 \\
\midrule
\multicolumn{6}{l}{\textbf{Variant B - Velocity and Position-Supervised Refinement}} \\
$\parallel \mathbf{r}_2 - \mathbf{r}_f \parallel$ (km) & 0.234 & 0.0635 & 0.12 & 0.212 & 1.7 \\
$\parallel \mathbf{v}_{1,true} - \mathbf{v}_{1}^{(K)} \parallel$ (m/s) & 3.05 & 0.0448 & 0.0873 & 0.179 & 25.9 \\
\midrule
\multicolumn{6}{l}{\textbf{Variant C - Position-Supervised Refinement}} \\
$\parallel \mathbf{r}_2 - \mathbf{r}_f \parallel$ (km) & 0.0645 & 0.0162 & 0.0273 & 0.0466 & 0.517 \\
$\parallel \mathbf{v}_{1,true} - \mathbf{v}_{1}^{(K)} \parallel$ (m/s) & 2.52 & 0.0111 & 0.0197 & 0.0374 & 9.27 \\
\bottomrule
\end{tabular}
\end{table}

\subsubsection{Multi Revolution LEO Results}

\begin{figure}[tbp]
    \centering
    \includegraphics[width=0.9\linewidth]{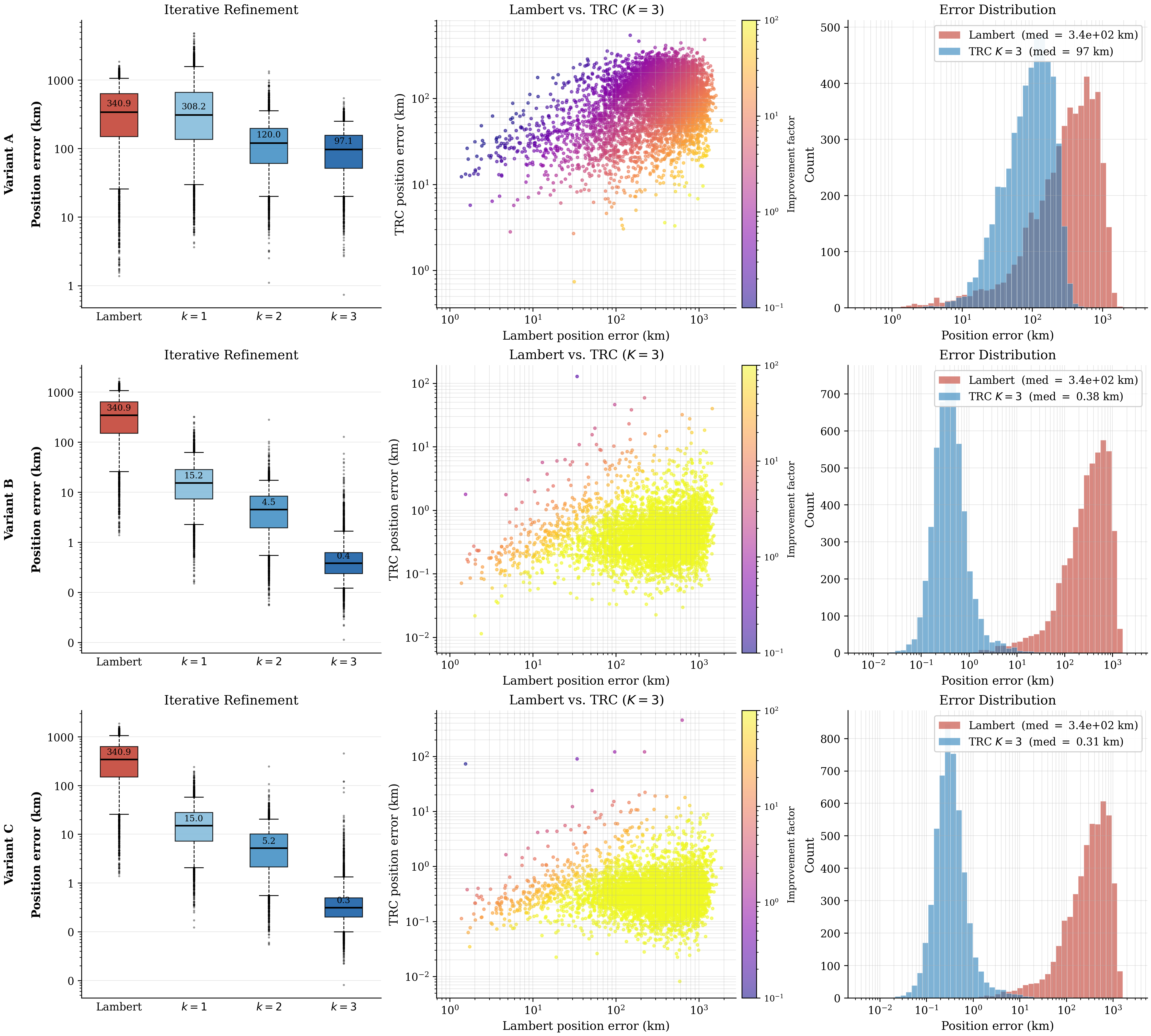}
    \caption{Performance of TRM on the multi-revolution LEO dataset}
    \label{fig:s2_multirev}
\end{figure}

The three variants are evaluated on the multi-revolution LEO validation set of $N = 5$k $J_2$-perturbed Lambert transfers spanning revolution counts up to $n_\mathrm{rev} = 13$. Results are reported in Table~\ref{tab:err_multi_rev} and Figure~\ref{fig:s2_multirev}.

Variants~B and~C both converge to sub-kilometre median accuracy. The position-only variant (C) reaches a median terminal-position error of 0.31\,km with a 99th percentile of 5.14\,km; the velocity-supervised variant (B) is comparable at 0.38\,km median and 5.81\,km at the 99th percentile. The end-to-end learned-Lambert variant (Variant~A) reaches only 97\,km median terminal error with a mean of 111\,km and a 99th percentile above 300\,km, placing the bulk of the distribution at the same order of magnitude as the uncorrected Lambert baseline. The performance gap between Variant~A and the refinement variants widens from two orders of magnitude on single-revolution LEO to nearly three orders here, consistent with the harder optimisation landscape introduced when the network must learn the multi-revolution Lambert solution itself in addition to the $J_2$ correction.

\begin{table}[tbp]
\centering
\small
\setlength{\tabcolsep}{5pt}
\caption{Final position and corrected initial velocity error distributions on the validation set for multi-revolution LEO. $Q_2$ is the median; $P_{99}$ is the 99th percentile.}
\label{tab:err_multi_rev}
\begin{tabular}{lrrrrr}
\toprule
Quantity & Mean & $Q_1$ & $Q_2$ & $Q_3$ & $P_{99}$ \\
\midrule
\multicolumn{6}{l}{\textbf{Variant A - Learned Lambert}} \\
$\parallel \mathbf{r}_2 - \mathbf{r}_f \parallel$ (km) & 111 & 52 & 97.1 & 157 & 316 \\
$\parallel \mathbf{v}_{1,true} - \mathbf{v}_{1}^{(K)} \parallel$ (m/s) & 98.4 & 25.4 & 64.9 & 141 & 444 \\
\midrule
\multicolumn{6}{l}{\textbf{Variant B - Velocity and Position-Supervised Refinement}} \\
$\parallel \mathbf{r}_2 - \mathbf{r}_f \parallel$ (km) & 0.703 & 0.239 & 0.384 & 0.628 & 5.81 \\
$\parallel \mathbf{v}_{1,true} - \mathbf{v}_{1}^{(K)} \parallel$ (m/s) & 2.44 & 0.165 & 0.293 & 0.594 & 20.5 \\
\midrule
\multicolumn{6}{l}{\textbf{Variant C - Position-Supervised Refinement}} \\
$\parallel \mathbf{r}_2 - \mathbf{r}_f \parallel$ (km) & 0.702 & 0.201 & 0.311 & 0.497 & 5.14 \\
$\parallel \mathbf{v}_{1,true} - \mathbf{v}_{1}^{(K)} \parallel$ (m/s) & 2.57 & 0.143 & 0.249 & 0.471 & 19.8 \\
\bottomrule
\end{tabular}
\end{table}

\subsubsection{Multi-Revolution Jovian Results}

\begin{figure}[tbp]
    \centering
    \includegraphics[width=\linewidth]{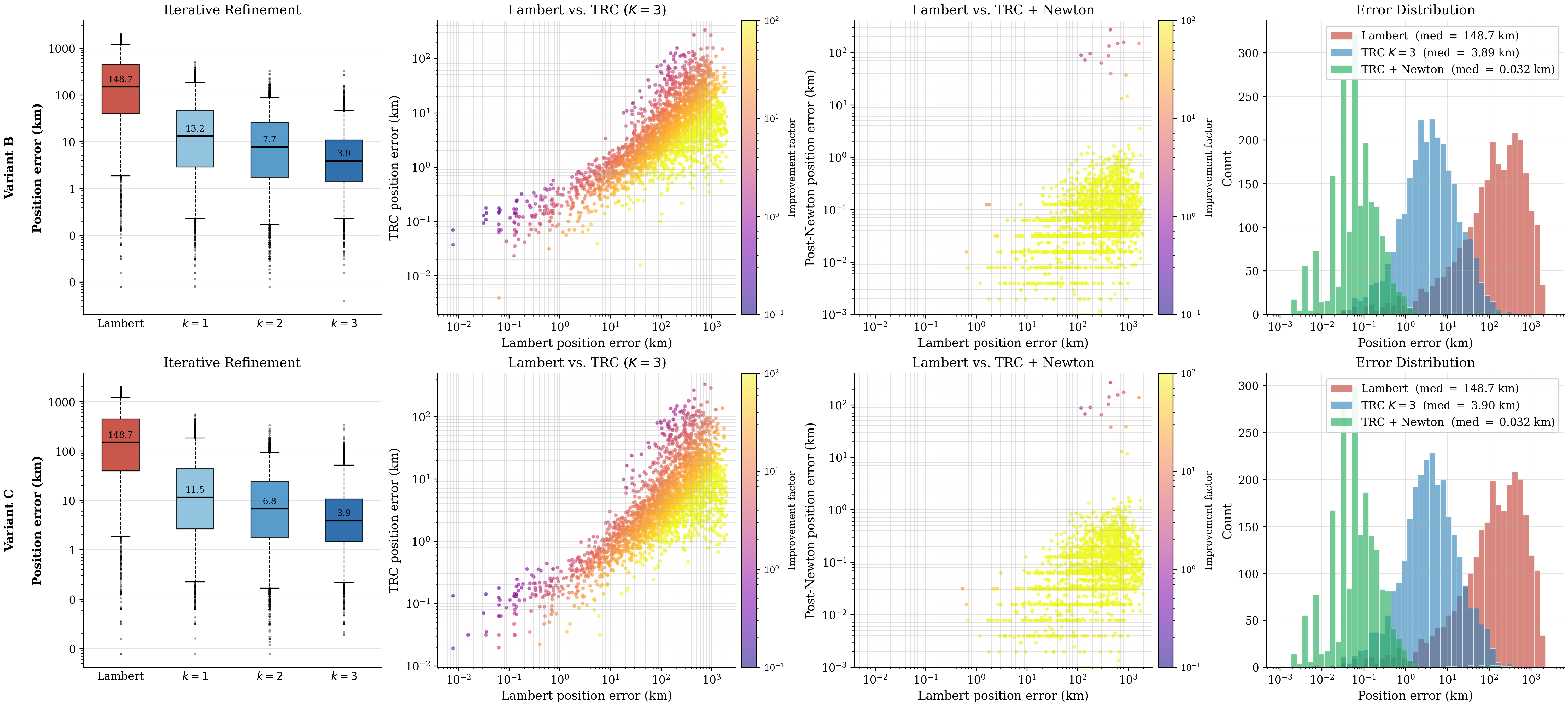}
    \caption{Performance of TRM on the multi-revolution Jovian dataset. Results also include the outcomes after a single Newton shooting iteration from the TRM output. }
    \label{fig:s2_jov}
\end{figure}

The Jovian validation set consists of $N = 20$k transfers from the dataset of Yang et al.~\cite{yang2022}, with revolution counts up to $n_\mathrm{rev} = 11$. This is the most difficult of the three test cases: time-of-flight ranges over $\sim 5\times 10^6$\,s (roughly 60 days), the $J_2$ correction is heavy-tailed, and the uncorrected Lambert baseline median terminal error is approximately 4500\,km which is two orders of magnitude larger than in either LEO case. Results are reported in Table~\ref{tab:err_jovian_thr10} and Figure~\ref{fig:s2_jov}.

Variant~A was not attempted for this case as it fails for the simpler LEO cases.  Variants~B and~C both reduce the median terminal-position error by three orders of magnitude relative to the uncorrected Lambert baseline, to 3.89\,km (Variant~B) and 3.90\,km (Variant~C). The two variants are statistically indistinguishable on the central distribution, with $Q_1 \approx 1.45$\,km, $Q_3 \approx 10.8$\,km, and means of $\sim 11$\,km. The means substantially exceed the medians because of a heavy upper tail, driven by chaotic high-revolution geometries near low perijove where the perturbed dynamics are most sensitive to initial-velocity uncertainty. Velocity-error statistics show the same pattern: medians of $\sim 0.04$\,m/s with 99th-percentile tails extending to 15.7\,m/s (B) and 26.2\,m/s (C).

Table~\ref{tab:err_jovian_thr10} also reports the effect of a single Newton-shooting
iteration applied on top of the TRM-PL velocity output, enabling a direct comparison with
the DNN-plus-Newton pipeline of Yang et al.~\cite{yang2022}. With one corrector step, the
median terminal error drops from $\sim 4$\,km to 0.063\,km on both variants, the
99th-percentile tail tightens from $\sim 100$\,km to $\sim 20$\,km, and the median velocity
error falls below $10^{-3}$\,m/s. By comparison, the DNN-only initial guess reported in
Yang et al.~\cite[Fig.~8]{yang2022} has a per-component terminal-position standard deviation
on the order of tens of kilometers, and their pipeline requires several Newton iterations to
reach the stated 1\,m tolerance, with the iteration count growing with revolution number.
The TRM-PL output is therefore a sufficiently accurate initial guess that a single Newton
iteration reaches sub-100-meter accuracy across the full validation distribution. This
comparison is necessarily indirect: Yang et al.~\cite{yang2022} report no distribution
statistics for their training and validation data, thus it cannot be guaranteed that the datasets, generated here using the procedure described in their paper, match theirs exactly.

\begin{table}[tbp]
\centering
\small
\setlength{\tabcolsep}{5pt}
\caption{Final position and corrected initial velocity error distributions on the validation set. $Q_2$ is the median; $P_{99}$ is the 99th percentile for the multi-revolution Jovian case. Each variant reports both the TRM-PL output and the result after 1 Newton refinement iteration on top of the TRM-PL velocity guess.}
\label{tab:err_jovian_thr10}
\begin{tabular}{lrrrrr}
\toprule
Quantity & Mean & $Q_1$ & $Q_2$ & $Q_3$ & $P_{99}$ \\
\midrule
\multicolumn{6}{l}{\textbf{Variant B - Velocity and Position-Supervised Refinement}} \\
$\parallel \mathbf{r}_2 - \mathbf{r}_f \parallel$ (km) & 10.7 & 1.43 & 3.89 & 10.8 & 95.2 \\
$\parallel \mathbf{v}_{1,true} - \mathbf{v}_{1}^{(K)} \parallel$ (m/s) & 0.742 & 0.0164 & 0.0428 & 0.136 & 15.7 \\
\addlinespace[2pt]
\multicolumn{6}{l}{\textit{After 1 Newton iteration}} \\
$\parallel \mathbf{r}_2 - \mathbf{r}_f \parallel$ (km) & 1.07 & 0.00781 & 0.0625 & 0.147 & 15.6 \\
$\parallel \mathbf{v}_{1,true} - \mathbf{v}_{1}^{(K)} \parallel$ (m/s) & 0.232 & 4.77e-04 & 9.83e-04 & 0.00226 & 0.197 \\
\midrule
\multicolumn{6}{l}{\textbf{Variant C - Position-Supervised Refinement}} \\
$\parallel \mathbf{r}_2 - \mathbf{r}_f \parallel$ (km) & 11.5 & 1.47 & 3.9 & 10.7 & 109 \\

$\parallel \mathbf{v}_{1,true} - \mathbf{v}_{1}^{(K)} \parallel$ (m/s) & 0.972 & 0.0165 & 0.0463 & 0.15 & 26.2 \\
\addlinespace[2pt]
\multicolumn{6}{l}{\textit{After 1 Newton iteration}} \\
$\parallel \mathbf{r}_2 - \mathbf{r}_f \parallel$ (km) & 1.22 & 0.00781 & 0.0625 & 0.156 & 25.1 \\
$\parallel \mathbf{v}_{1,true} - \mathbf{v}_{1}^{(K)} \parallel$ (m/s) & 0.238 & 4.77e-04 & 9.83e-04 & 0.00227 & 0.232 \\
\bottomrule
\end{tabular}
\end{table}

\section{Conclusions}
\label{sec:conclusion}

This paper presents TRM-PL, a recursive neural solver for the $J_2$-perturbed Lambert problem based on the TRM architecture.
 The solver applies a shared two-level reasoning module repeatedly inside a single end-to-end differentiable network, where at each iteration, the current departure-velocity estimate is propagated under $J_2$ through a differentiable RK4 integrator, the resulting terminal-position residual is fed back into the network, and a residual decoder produces a velocity correction.  Three training paradigms were compared on the same architecture: (A) jointly learning the Keplerian Lambert solution and the $J_2$ correction; (B) refining from the classical Lambert solution with both position and true perturbed initial velocity supervision; and (C) refining from the classical Lambert solution with position-only supervision. The paradigms were evaluated on three test cases of increasing difficulty: single-revolution LEO transfers, multi-revolution LEO transfers up to thirteen revolutions, and multi-revolution Jovian transfers up to eleven revolutions based on the dataset of Yang et al.~\cite{yang2022}.

The findings highlight that the position-only supervision variant (C) is the most accurate and stable across all three regimes, reducing the median terminal-position error of the uncorrected Lambert baseline by two to three orders of magnitude. The joint-training variant (A) degrades sharply with problem difficulty and fails to converge on the Jovian dataset, indicating that asking the network to also reproduce the Keplerian Lambert solution is a poor inductive bias when an exact, inexpensive solver is available; the recursive network's strength is in iterative refinement, not in approximating the Lambert solution itself. The recursive refinement loop can be seen as a learned alternative to the homotopy and continuation schemes used in classical perturbed-Lambert solvers: each iteration is one step of a continuation path discovered from data rather than prescribed in advance. The training of this loop is reinforcement-learning-like, in that the network is rewarded for reducing the terminal-position error rather than supervised against a fixed velocity label. TRM-PL provides a unified, end-to-end differentiable replacement for the two-stage DNN-plus-Newton-shooting pipeline of Yang et al.~\cite{yang2022} for the perturbed Lambert Problem. On their multi-revolution Jovian dataset, TRM-PL reaches a median terminal-position error of $\sim 4$\,km in a single forward pass, and drops to  0.063\,km after a single Newtown corrector iteration,  performing an order of magnitude better.


\bibliographystyle{unsrtnat}
\bibliography{references}

 \appendix
\section*{Appendix: Training and Validation Dataset Statistics}

\begin{table}[tbp]
\centering
\small
\caption{Sample statistics for the three training datasets, including mean, standard deviation, 10th percentile $p_{10}$ and 90th percentile $p_{90}$}
\label{tab:stats_all}
\small
\begin{tabular}{lcccc}
\toprule
Quantity & Mean & Std & $p_{10}$ & $p_{90}$ \\
\midrule
\multicolumn{5}{c}{\textbf{LEO Single Revolution} (50k)} \\
$\Delta t$ (hr) 
& $0.870$ 
& $0.411$ 
& $0.325$ 
& $1.414$ \\

$\mathbf{v}_\mathrm{1, Lambert}$ (km/s) 
& $\begin{bmatrix}-0.02564\\ -0.002255\\ 0.009541\end{bmatrix}$ 
& $\begin{bmatrix}4.995\\ 3.544\\ 3.547\end{bmatrix}$ 
& $\begin{bmatrix}-6.753\\ -4.952\\ -4.944\end{bmatrix}$ 
& $\begin{bmatrix}6.725\\ 4.951\\ 4.966\end{bmatrix}$ \\

$\mathbf{v}_{1, \mathrm{true}}$ (km/s) 
& $\begin{bmatrix}-0.02553\\ -0.001734\\ 0.0098\end{bmatrix}$ 
& $\begin{bmatrix}4.995\\ 3.546\\ 3.547\end{bmatrix}$ 
& $\begin{bmatrix}-6.752\\ -4.956\\ -4.944\end{bmatrix}$ 
& $\begin{bmatrix}6.723\\ 4.957\\ 4.963\end{bmatrix}$ \\

$\mathbf{e}_{\mathrm{Lambert}}$ (km) 
& $\begin{bmatrix}0.008508\\ -0.01862\\ -0.1065\end{bmatrix}$ 
& $\begin{bmatrix}33.33\\ 25.41\\ 23.89\end{bmatrix}$ 
& $\begin{bmatrix}-37.44\\ -23.79\\ -22.67\end{bmatrix}$ 
& $\begin{bmatrix}37.36\\ 23.57\\ 22.78\end{bmatrix}$ \\
\midrule
\multicolumn{5}{c}{\textbf{LEO Multi Revolution} (50k)} \\
$\Delta t$ (hr) 
& $11.006$ 
& $5.628$ 
& $2.864$ 
& $18.325$ \\
$\mathbf{v}_\mathrm{1,Lambert}$ (km/s) 
& $\begin{bmatrix}0.003572\\ -0.0135\\ 0.01248\end{bmatrix}$ 
& $\begin{bmatrix}5.135\\ 3.606\\ 3.602\end{bmatrix}$ 
& $\begin{bmatrix}-6.79\\ -5.099\\ -5.093\end{bmatrix}$ 
& $\begin{bmatrix}6.806\\ 5.092\\ 5.113\end{bmatrix}$ \\

$\mathbf{v}_{1, \mathrm{true}}$ (km/s) 
& $\begin{bmatrix}0.003609\\ -0.01322\\ 0.01195\end{bmatrix}$ 
& $\begin{bmatrix}5.135\\ 3.616\\ 3.595\end{bmatrix}$ 
& $\begin{bmatrix}-6.789\\ -5.116\\ -5.083\end{bmatrix}$ 
& $\begin{bmatrix}6.805\\ 5.105\\ 5.103\end{bmatrix}$ \\
$\mathbf{e}_{\mathrm{Lambert}}$ (km) 
& $\begin{bmatrix}0.9106\\ 0.3511\\ 0.2508\end{bmatrix}$ 
& $\begin{bmatrix}367.8\\ 286.7\\ 255.5\end{bmatrix}$ 
& $\begin{bmatrix}-450.5\\ -312.1\\ -278.7\end{bmatrix}$ 
& $\begin{bmatrix}452.3\\ 311.3\\ 280.8\end{bmatrix}$ \\
\midrule
\multicolumn{5}{c}{\textbf{Jovian Multi Revolution} (200k)} \\
$\Delta t$ (hr) 
& $1431.389$ 
& $1056.111$ 
& $224.528$ 
& $3000.000$ \\

$\mathbf{v}_{1,\mathrm{Lambert}}$ (km/s) 
& $\begin{bmatrix}0.03182\\ -0.0219\\ 0.00386\end{bmatrix}$ 
& $\begin{bmatrix}6.431\\ 6.413\\ 3.617\end{bmatrix}$ 
& $\begin{bmatrix}-7.96\\ -7.988\\ -4.717\end{bmatrix}$ 
& $\begin{bmatrix}8.018\\ 7.966\\ 4.723\end{bmatrix}$ \\

$\mathbf{v}_{1,\mathrm{true}}$ (km/s) 
& $\begin{bmatrix}0.03182\\ -0.0218\\ 0.003942\end{bmatrix}$ 
& $\begin{bmatrix}6.431\\ 6.414\\ 3.617\end{bmatrix}$ 
& $\begin{bmatrix}-7.96\\ -7.988\\ -4.72\end{bmatrix}$ 
& $\begin{bmatrix}8.018\\ 7.966\\ 4.726\end{bmatrix}$ \\
$\mathbf{e}_{\mathrm{Lambert}}$ (km) 
& $\begin{bmatrix}-8.461\\ 15.46\\ 3.936\end{bmatrix}$ 
& $\begin{bmatrix}4648\\ 4581\\ 2604\end{bmatrix}$ 
& $\begin{bmatrix}-3710\\ -3729\\ -2275\end{bmatrix}$ 
& $\begin{bmatrix}3729\\ 3734\\ 2276\end{bmatrix}$ \\

\bottomrule
\end{tabular}
\end{table}

\begin{table}[tbp]
\centering
\small
\caption{Sample statistics for the validation dataset, including mean, standard deviation, 10th percentile $p_{10}$, and 90th percentile $p_{90}$.}
\label{tab:stats_all_val}
\begin{tabular}{lcccc}
\toprule
Quantity & Mean & Std & $p_{10}$ & $p_{90}$ \\
\midrule

\multicolumn{5}{c}{\textbf{LEO Single Revolution} (5k)} \\
$\Delta t$ (hr) 
& $0.873$ 
& $0.411$ 
& $0.324$ 
& $1.414$ \\

$\mathbf{v}_{1,\mathrm{Lambert}}$ (km/s) 
& $\begin{bmatrix}0.0088\\ -0.002156\\ 0.03086\end{bmatrix}$ 
& $\begin{bmatrix}4.997\\ 3.512\\ 3.585\end{bmatrix}$ 
& $\begin{bmatrix}-6.685\\ -4.824\\ -4.892\end{bmatrix}$ 
& $\begin{bmatrix}6.684\\ 4.961\\ 4.996\end{bmatrix}$ \\

$\mathbf{v}_{1,\mathrm{true}}$ (km/s) 
& $\begin{bmatrix}0.008659\\ -0.003748\\ 0.03081\end{bmatrix}$ 
& $\begin{bmatrix}4.997\\ 3.515\\ 3.584\end{bmatrix}$ 
& $\begin{bmatrix}-6.684\\ -4.824\\ -4.901\end{bmatrix}$ 
& $\begin{bmatrix}6.681\\ 4.967\\ 4.996\end{bmatrix}$ \\

$\mathbf{e}_{\mathrm{Lambert}}$ (km) 
& $\begin{bmatrix}0.2166\\ 0.08838\\ 0.14\end{bmatrix}$ 
& $\begin{bmatrix}33.54\\ 25.63\\ 25.82\end{bmatrix}$ 
& $\begin{bmatrix}-39.88\\ -23.21\\ -22.57\end{bmatrix}$ 
& $\begin{bmatrix}38.43\\ 22.28\\ 22.27\end{bmatrix}$ \\

\midrule
\multicolumn{5}{c}{\textbf{LEO Multi Revolution} (5k)} \\
$\Delta t$ (hr) 
& $11.042$ 
& $5.608$ 
& $2.770$ 
& $18.236$ \\

$\mathbf{v}_{1,\mathrm{Lambert}}$ (km/s) 
& $\begin{bmatrix}0.04084\\ 0.08484\\ -0.05466\end{bmatrix}$ 
& $\begin{bmatrix}5.135\\ 3.601\\ 3.619\end{bmatrix}$ 
& $\begin{bmatrix}-6.798\\ -5.061\\ -5.116\end{bmatrix}$ 
& $\begin{bmatrix}6.768\\ 5.057\\ 5.101\end{bmatrix}$ \\

$\mathbf{v}_{1,\mathrm{true}}$ (km/s) 
& $\begin{bmatrix}0.0406\\ 0.08304\\ -0.05422\end{bmatrix}$ 
& $\begin{bmatrix}5.135\\ 3.611\\ 3.612\end{bmatrix}$ 
& $\begin{bmatrix}-6.799\\ -5.045\\ -5.145\end{bmatrix}$ 
& $\begin{bmatrix}6.773\\ 5.07\\ 5.102\end{bmatrix}$ \\
$\mathbf{e}_{\mathrm{Lambert}}$ (km) 
& $\begin{bmatrix}-2.172\\ 7.012\\ -1.693\end{bmatrix}$ 
& $\begin{bmatrix}369.6\\ 288.6\\ 252.5\end{bmatrix}$ 
& $\begin{bmatrix}-452.3\\ -311.9\\ -283\end{bmatrix}$ 
& $\begin{bmatrix}459.8\\ 322.2\\ 275.9\end{bmatrix}$ \\
\bottomrule
\multicolumn{5}{c}{\textbf{Jovian Multi Revolution } (20k)} \\
$\Delta t$ (hr) 
& $1440.833$ 
& $1068.056$ 
& $217.111$ 
& $3038.889$ \\

$\mathbf{v}_{1,\mathrm{Lambert}}$ (km/s) 
& $\begin{bmatrix}-0.08684\\ -0.03257\\ -0.0349\end{bmatrix}$ 
& $\begin{bmatrix}6.355\\ 6.447\\ 3.622\end{bmatrix}$ 
& $\begin{bmatrix}-7.98\\ -7.98\\ -4.772\end{bmatrix}$ 
& $\begin{bmatrix}7.914\\ 7.991\\ 4.683\end{bmatrix}$ \\

$\mathbf{v}_{1,\mathrm{true}}$ (km/s) 
& $\begin{bmatrix}-0.08647\\ -0.03313\\ -0.03493\end{bmatrix}$ 
& $\begin{bmatrix}6.355\\ 6.447\\ 3.623\end{bmatrix}$ 
& $\begin{bmatrix}-7.98\\ -7.977\\ -4.777\end{bmatrix}$ 
& $\begin{bmatrix}7.914\\ 7.995\\ 4.68\end{bmatrix}$ \\
$\mathbf{e}_{\mathrm{Lambert}}$ (km) 
& $\begin{bmatrix}-2.877\\ -38.56\\ -3.622\end{bmatrix}$ 
& $\begin{bmatrix}4567\\ 4309\\ 2609\end{bmatrix}$ 
& $\begin{bmatrix}-3742\\ -3664\\ -2284\end{bmatrix}$ 
& $\begin{bmatrix}3676\\ 3660\\ 2274\end{bmatrix}$ \\

\bottomrule

\end{tabular}
\end{table}

\end{document}